\documentclass{article}
\usepackage[T2A]{fontenc}
\usepackage[cp1251]{inputenc}
\usepackage[english]{babel}

\begin{document}

\begin{center} {\bf \Large
Probabilistic counterparts of  nonlinear parabolic PDE systems.}\end{center} 
\begin{center}Ya.I. Belopolskaya, \\
St.Petersburg State University for
Architecture and Civil Engineering,\\
St.Petersburg, Russia, 
yana@yb1569.spb.edu \end{center} 


\abstract{We extend the results of the FBSDE theory in order  to
construct a probabilistic representation of a viscosity solution to
the Cauchy problem for a system  of quasilinear parabolic equations.
We  derive a BSDE associated with a class of quailinear parabolic
system and prove the existence and uniqueness of its solution.  To
be able to deal with systems including nondiagonal first order terms
along with the underlying diffusion process we consider its
multiplicative operator functional. We essentially exploit as well
the fact that the system under consideration can be reduced to a
scalar equation in a enlarged phase space. This allows to obtain
some comparison theorems and to prove that a solution to FBSDE gives
rise to a viscosity solution of the original Cauchy   problem for a
system of quasilinear parabolic equations. }

\section{Introduction}
\label{sec:1}

Quasilinear systems of parabolic equations   arise as mathematical models which  describe various chemical and  biological phenomena. They arise  as well in financial mathematics and in differential geometry when one considers nonlinear parabolic equations in sections of vector bundles.

Let $d,d_1$ be given integers, $a(x)\in R^d, A(x)\in R^{d\times d},
B(x)\in R^{d\times d_1\times d_1}, c(x)\in R^{d_1\times d_1}, x\in
R^d$ and $g:R^d\times R^{d_1}\times R^{d\times d_1}\to {R^d_1}$ be
given. Consider a class of quasilinear parabolic equations of the
form
\begin{equation}\label{1.1}
\frac{\partial u_l}{\partial s}+\frac 12 Tr A^*\nabla^2 u_lA+\langle
a,\nabla u_l\rangle+
\end{equation}
$$
+B_{lm}^i\nabla_i u_m+c_{lm}u_m +g_l(s,x,u,\nabla u)=0,\quad
u_l(T,x)=u_{0l}(x),\quad l=1,\dots,d_1
$$
with respect of $R^{d_1}$-valued function $u(s,x)$ defined on
$[0,T]\times R^d.$ Here and below we assume a convention of summing
up over repeating indices if the contrary is not mentioned and
denote by  $\langle\cdot,\cdot\rangle$   an inner product in $R^d$
regardless of $d$.

 One can suggest at least a couple of  probabilistic counterparts of the Cauchy problem (\ref{1.1}). To derive them let us
  assume first that there exists a classical  solution $u(s,x)$ to this problem.
 In this case one can  prove applying the standard technique of the stochastic differential equation theory and especially the Ito formula, that the function $u(s,x)$ satisfying (\ref{1.1}) admits at least two probabilistic representations.

 The first one was suggested in papers by Dalecky and Belopolskaya [1] -[3] and was aimed to  develop a probabilistic
 approach to prove the existence and uniqueness of a classical solution to (\ref{1.1}) and as well as to much more general systems of the form
$$\frac{\partial u_l}{\partial s} +F(x,u,\nabla u,\nabla ^2u_l)=0,\quad u_l(T,x)=u_{l0}(x).$$

 The second one  suggested in papers by Pardoux and Peng  [4]- [6]  leads to the powerful backward
  stochastic differential equations (BSDE) theory. This approach allows to construct a viscosity solution to
  a quasilinear scalar parabolic PDE or to  a diagonal system of PDEs (see [6] -- [7]).
   In terms of (\ref{1.1}) this  means that one have to set  $B\equiv 0$ and $c\equiv 0$ and
   $g_l(x,u, A^*\nabla u)\equiv g_l(x,u, A^*\nabla u_l)$.

 To present these approaches we fix
 a probability space $(\Omega,{\cal F},P)$ and  denote by $w(t)\in R^d$ the standard Wiener process. Let ${\cal F}_t$ be a flow of $\sigma$-subalgebras of ${\cal F}$ generated by $w(t)$ and   $E_{s,x} [f(\xi(T))]=E[f(\xi(T)|\xi(s)=x]$ denote the conditional expectation.

 Assume that  $g$ in (\ref{1.1}) does not depend on $\nabla u$  and all coefficients $a,A,B,C$ depend on $s,x$ and $u$. Assume that   $u(s,x)$ is  a   smooth  function  satisfying (\ref{1.1}) with these parameters. Then it  was stated in [1]  that this function
 admits a representation of the form
 \begin{equation}\label{1.2}
\langle h,u(s,x)\rangle=E_{s,x}\left[\langle
\eta(T),u_0(\xi(T))\rangle +\int_s^T\langle
\eta(\theta),g(\theta,\xi(\theta), u(\theta,
\xi(\theta)))d\theta\right],
\end{equation}
where  stochastic processes    $\xi(t)$  and $\eta(t)$ satisfy the stochastic equations
 \begin{equation}\label{1.3}
d\xi(t)=a(\xi(t),u(t,\xi(t)))dt +A(\xi(t),u(t,\xi(t)))dw(t),\quad
\xi(s)=x,
\end{equation}
and
 \begin{equation}\label{1.4}
d\eta(t)=c(\xi(t),u(t,\xi(t)))\eta(t)dt
+C(\xi(t),u(t,\xi(t)))(\eta(t),dw(t)),\quad \eta(s)=h.
\end{equation}
Note that  $a,A,c$ in (\ref{1.3}), (\ref{1.4}) are the same as in (\ref{1.1}) while  it is assumed that $C$ in (\ref{1.4}) and $B$  in (\ref{1.1})
satisfy an equality   $B^{lm}_k =C^{lm}_iA_{ik}.$

{\bf Remark.}
 Notice that when $A$ is a nondegenerated matrix one can  define $C$ by
$C^{lm}_i=B^{lm}_kA^{-1}_{ki}$ while when $A$ is a degenerated matrix we assume that $B$ has the above form.
It is important that although  in a general case  $A=0$ yields  $B=0$ and we do not obtain a general  nonlinear system of hyperbolic equations as a vanishing viscosity limit of (\ref{1.1}). Nevertheless one can state some restrictions on
$B$  such that   given $A_{\epsilon}=\epsilon A $ and $C_{\epsilon}=\epsilon^{-1}C$  one can apply (\ref{1.2}) to investigate the vanishing viscosity limit of (\ref{1.1}) with these coefficients  (see {AB}) .

  An important observation is the fact that we can  consider (\ref{1.2})-(\ref{1.4}) as a closed system of equations and
   state conditions on its data to ensure the existence and uniqueness of a solution to this  system. If in addition it will be revealed that the function $u(s,x)$ given by  (\ref{1.2}) is twice differentiable in the spatial variable $x$, then one can verify that $u(s,x)$ is a unique classical solution of  (\ref{1.1}) with correspondent parameters.  It should be mentioned that this approach can be essentially generated to give a possibility to study systems of quasilinear and even fully nonlinear parabolic equations. In other words one can consider  (\ref{1.1})  with  coefficients $ a,A,c,C, g$  depending on  $(x,u,\nabla u)$ or even $(x,u,\nabla u,\nabla^2 u)$. Note that to deal with these more complicated cases within a framework of this approach we require more strong assumptions
concerning regularity of coefficients of(\ref{1.3})-(\ref{1.4}) and
the Cauchy data $u_0$. As a result we can  prove  on this way the
existence and uniqueness of a classical solution to (\ref{1.1}),
possibly on a small time interval.

To describe the second  approach which allows to construct a different class of solutions to the Cauchy problem
 \begin{equation}\label{1.5}
 \frac{\partial u_l}{\partial s}+\frac 12 Tr A^*(x)\nabla^2 u_lA(x)+\langle a(x),\nabla u_l\rangle
+g_l(x,u, A^*\nabla u_l)=0,\quad u_l(T,x)=u_{0l}(x),
\end{equation}
 we assume once again that there exists a classical solution $u_l(s,x)$ of  (\ref{1.5}).

 Consider a stochastic process $\xi(t)$ satisfying  (\ref{1.3})  with coefficients
 $a(s,x,u)\equiv a(s,x),$ $ A(s,x,u)\equiv A(s,x)$. Keeping in mind that $u_l(s,x)$ is a classical solution of
 (\ref{1.5}), by Ito's formula  we derive an    expression for a stochastic differential of $y(t)=\Gamma^*(s,t)u(t,\xi(t))$ in the form
 \begin{equation}\label{1.6}
dy(t)=-g(t,\xi(t),y(t),z(t))dt-zdw(t),\quad y(T)=\Gamma^*(s,T)u_0(\xi(T)),
\end{equation}
where $z(t)=A^*(\xi(t))\nabla u(t,\xi(t)), \eta(t)=\Gamma(s,t)h.$ The equation (\ref{1.6})
is called a backward stochastic equation (BSDE).

In general one can forget about the process $\xi(t)$ and consider an independent BSDE of the form
  \begin{equation}\label{1.7}
dy(t)=-f(t,y(t),z(t))dt-zdw(t),\quad y(T)=\zeta,
\end{equation}
where $f(t,y,z)$ is an ${\cal F}_t$-adapted random process meeting some additional requirements
and $\zeta$ is an  ${\cal F}_T$-measurable random variable.  A general theory of BSDEs was developed
by a number of authors (see  e.g. [7]   for references).  In addition the system  (\ref{1.4}), (\ref{1.6})
 shows a way to construct
  the so called viscosity solution   to (\ref{1.5}) (defined in [8])  setting  $u(s,x)=y(s))$.

To generalize this approach and apply it to
(\ref{1.1})  we observe that this system has a
crucial property which can be easily revealed if one analyzes the probabilistic representation (\ref{1.2}) of a smooth solution to (\ref{1.1}). Namely, the Cauchy problem  (\ref{1.1})
can be reduced to the  Cauchy problem  for a scalar equation
 \begin{equation}\label{1.8}
\frac{\partial \Phi}{\partial s}+\frac 12 Tr Q^*(x,h)\nabla^2 \Phi
Q(x,h)+\langle q(x,h),\nabla \Phi\rangle +G(s,h, x,\Phi,
Q^*\nabla\Phi)=0,\quad
\end{equation}
$$\Phi(T,x)=\Phi_0(x,h)=\langle h,u_{0}(x)\rangle.$$
with respect to a scalar function $\Phi(s,x,h)=\langle h,u(s,x)\rangle.$

Here
$$Tr Q^*\nabla^2\Phi(s,x,h)Q= A^*_{ki}\frac{\partial^2\Phi(s,x,h)}{\partial x_i\partial x_j}A_{jk}+
2C^{lm}_kh_l\frac{\partial^2\Phi(s,x,h)}{\partial x_j\partial
h_m}A_{jk}+ $$ $$+C^{qm}_kh_m\frac{\partial^2\Phi(s,x,h)}{\partial
h_q\partial h_p}C^{pn}_kh_n=
A^*_{ki}\frac{\partial^2\Phi(s,x,h)}{\partial x_i\partial
x_j}A_{jk}+ 2C^{lm}_kh_l\frac{\partial^2\Phi(s,x,h)}{\partial
x_j\partial h_m}A_{jk},$$ since, due to linearity of $\Phi(s,x,h)$
in $h$,  we have $\frac{\partial^2\Phi(s,x,h)}{\partial h_q\partial
h_p}\equiv 0.$ In addition $$ \langle q,\nabla\Phi(s,x,h)\rangle
=a_j\frac{\partial \Phi(s,x,h)}{\partial x_j}+ c_{lm}h_m\frac
{\partial \Phi(s,x,h)}{\partial h_l}, \quad G(s,x,h)=\langle
h,g(s,x,u, A^*\nabla u)\rangle.$$

Coming back to (\ref{1.4})   we notice that its solution (provided it exists)
gives rise to a multiplicative operator functional  $\Gamma(t,s,\xi(\cdot))\equiv \Gamma (t,s)$ of the process  $\xi(t)$ satisfying (\ref{1.3}), that is $\eta(t)=\Gamma(t,s)h$ and $\Gamma(t,s)h=\Gamma(t,\theta)\Gamma(\theta,s)$  a.s. for $0\le s\le \theta\le t\le T$.
 Hence to derive an FBSDE  associated with   (\ref{1.1})   we can proceed as follows.

 Assume that there exists a classical solution to the Cauchy problem  (\ref{1.1})  or what is equivalent suppose that there exists a classical solution to    (\ref{1.8}) and compute  a stochastic differential of a stochastic process $Y(t)=\langle\eta(t), u(t,\xi(t))\rangle $,
$$dY(t)= \langle d\eta(t), u(t,\xi(t))\rangle+\langle\eta(t), du(t,\xi(t))\rangle+\langle d\eta(t),du(t,\xi(t))\rangle.$$
Taking into account  (\ref{1.3}),  (\ref{1.4}) by Ito's formula we derive the relation
 \begin{equation}\label{1.9}
dY(t)=- F(t,Y(t),Z(t))dt+ \langle Z(t),dW(t)\rangle,\quad Y(T)=\zeta= \langle \eta(T),u_0(\xi(T))\rangle,
\end{equation}
where $W(t)=(w(t),w(t))^*$,
$$\langle  Z(t), dW(t)\rangle=\langle C(\Gamma (t)h, dw(t)), u(t,\xi(t))\rangle +
\langle\Gamma(t)h,\nabla u(t,\xi(t))A dw\rangle=$$$$=\langle h,
\Gamma^*(t)[C^*u(t,\xi(t))+A^*\nabla u(t,\xi(t))]dw(t)\rangle$$ and
$\Gamma(t)h\equiv \Gamma(t,s)h=\eta_{s,h}(t).$ As a result we can
rewrite  (\ref{1.9})   in the form
 \begin{equation}\label{1.10}
dy(t) = -f(t,y(t), z(t))dt+ z(t) dw(t), \quad y(T)=\Gamma^*(s,T)u_0(\xi(T)),
\end{equation}
where
\begin{equation}\label{1.11}f(t,y(t),z(t))= \end{equation}$$=\Gamma^* (t)g\left(\xi(t),u(t,\xi(t)),C^*(t,\xi(t))u(t,\xi(t))+
  A^*(t,\xi(t))\nabla u(t,\xi(t))\right)=$$$$
= \Gamma^*
(t)g\left(\xi(t),[\Gamma^*]^{-1}(t)y(t),C^*(\xi(t))[\Gamma^*]^{-1}(t)y(t)+
 A^* (\xi(t))  [\Gamma^*]^{-1}(t)z(t)\right),$$
 $$Z(t)= ([\Gamma^*]^{-1}(t)C^*(t,\xi(t))u(t,\xi(t)), [\Gamma^*]^{-1}(t) A^*(\xi(t))\nabla u(t,\xi(t)))^*,$$
 $$ z(t) dw(t)=[\Gamma^*]^{-1}(t)[C^* udw(t)+A^*\nabla u dw(t)] \in R^{d_1} $$
 and $\langle h, z(t) dw(t)\rangle =\langle  Z(t), dW(t)\rangle.$

 When the solution  $y(t)$ is a scalar process and a comparison theorem holds  one can prove
that the function $u(s,x)$ defined by $y(s)=u(s,x)$ is a viscosity
solution of the Cauchy problem for a corresponding quasilinear
parabolic equation. In a multidimensional case it was shown in [9]
 that  given a solution of the BSDE
\begin{equation}\label{1.12}
dy_l(t) = -g_l(t,\xi(t),y(t), z_l(t))dt+ \langle z_l(t),
dw(t)\rangle, \quad y(T)= \Gamma^*(s,T)u_{0}(\xi(T)),
\end{equation}
where $\xi(t)$ satisfies (\ref{1.9}) under some condition one can
prove that the function $u(s,x)=y(s)$ is a viscosity solution to the
Cauchy problem
\begin{equation}\label{1.13}
\frac{\partial u_l}{\partial s}+\frac 12 Tr A^*\nabla^2 u_lA+\langle
a,\nabla u_l\rangle +g_l(s,x,u,A^*\nabla u_l)=0,\quad
u_l(T,x)=u_{0l}(x).
\end{equation}

In this paper we show that a  certain combination of  two approaches  allows to   extend the results
 of forward -backward stochastic equations (FBSDEs) theory to construct a viscosity solution to the system
 of the form (\ref{1.1}). In particular we define  the very notion of a viscosity solution for (\ref{1.1})
 and  prove a comparison theorem for solutions of multidimensional BSDEs which is a crucial point in
 construction of the viscosity solution via a solution to a BSDE.

In the next section we give a construction of an FBSDE required to construct a viscosity solution for  (\ref{1.1}),     assuming that coefficients $a,\sigma, C, c$ do not depend  on $u$.  We  state here  conditions on the BSDE  parameters that ensure the existence and uniqueness of its solution. In section 3 we prove  a comparison theorem and in section 4 we state the notion of a viscosity solution of the Cauchy problem for (\ref{1.1})   and prove that FBSDE solution gives rise to a viscosity solution for (\ref{1.1}).

\section{Forward-backward stochastic differential equations}
\label{sec:2}
\setcounter{equation}{0}
% Always give a unique label
% and use \ref{<label>} for cross-references
% and \cite{<label>} for bibliographic references
% use \sectionmark{}
% to alter or adjust the section heading in the running head
%\setcounter{equation}{0}

In this section we introduce  notations and present in a suitable form necessary results from  FBSDE theory adapted to the case under consideration.

Given  integers $d, d_1$ consider Euclidian spaces $R^d, R^{d_1}$ and let $\|\cdot\|$  denote a norm in  $R^d$ and $\langle\cdot,\cdot\rangle$ denote an inner product regardless of $d$.

  Given a Euclidian  space $X$ let

$ \bullet$ $L^{p}_{t} (X) $ be a set of ${\cal F}_{t}$ -measurable $X$-valued   random variables,   $ E\|\xi\|^{p}<\infty$;

$ \bullet$  ${\cal H}^2_c(X)$ be a set of ${\cal F}_{t}$ -measurable $X$-valued semimartingales  such that \\
$E\left[\sup_{0\le t\le T}\|y(t)\|^2\right]<\infty;$

$ \bullet$  ${\cal H}^2_{t}(X)$  be a set of ${\cal F}_{s,t}$ -measurable $X$-valued semimartingales  such that \\
$E\left[\sup_{0\le\theta\le t}\|y(\theta)\|^2\right]<\infty;$

$ \bullet$ ${\cal H}^{2}(X)$  be a set of square integrable progressively measurable processes
 $z(t)\in X$ such that $ \quad E\left[\int_0^T\|z(\tau)\|^2d\tau\right]<\infty;$

$ \bullet$ ${\cal S}^2={\cal H}^2_{c}(R^{d_1})\cup {\cal H}^{2}(R^{d\times d_1});$

$ \bullet$ ${\cal S}^3={\cal H}^2_{c}(R^{d})\cup {\cal H}^2_{c}(R^{d_1})\cup {\cal H}^{2}_{T}(R^{d\times d_1});$

$ \bullet$ ${\cal B}^2= {\cal H}^2(R^{d_1})\cup {\cal H}^{2}_{T}(R^{d\times d_1});$

$ \bullet$ ${\cal B}^3={\cal H}^2(R^{d})\cup {\cal H}^2(R^{d_1})\cup {\cal H}^{2}_{T}(R^{d\times d_1});$

$L(R^d)$ be the space of bounded linear maps acting in $R^d;$

$L(R^d;R^{d_1})\equiv R^{d\times d_1} $  be the space of bounded linear maps acting from $R^d$ to $R^{d_1}$;

$ \bullet$
Given  $\beta>0$   and  $\phi\in   {\cal H}^{2}_{T}(R^{d})$   let  $\|\phi\|_{\beta}^{2}
 =E\left[\int_{0}^{T}e^{\beta t}\|\phi(t)\|^{2}dt\right] $   and  ${\cal
 H}^{2}_{T,\beta}(R^{d})$     be the space   ${\cal H}^{2,d}_{T}$ equipped with the norm
 $\|\cdot\|_{\beta}.$

Let $W(t)=(w(t),w(t))\in R^d\times R^d$
and  $\kappa(t)=(\xi(t),\eta(t))\in R^d\times R^{d_1}$  be a solution of
a system of SDEs
 \begin{equation}\label{2.1}
d\xi(t)= a(t,\xi(t))dt+A(t,\xi(t))dw(t),\quad \xi(s)=x\in R^d,
\end{equation}
\begin{equation}\label{2.2}
d\eta(t)= c(t,\xi(t))\eta(t)dt+ C(t,\xi(t))(\eta(t), dw(t)) ,\quad \eta(s)=h\in R^{d_1}.
\end{equation}
We say that condition {\bf C 2.1} holds if
coefficients  $a:[0,\infty)\times R^d\to R^d,\quad A:[0,\infty)\times R^d\to  L(R^{d}),$
 $c:[0,\infty)\times R^d\to L(R^{d_1}), \quad C:[0,\infty)\times R^d\to L(R^d; L(R^{d_1}))$ are  continuous in $t\in [0,T]$ and  there exist constants $  K_1,K_2, L_1,L_2$ such that
$$\|a(t,x)\|^2+\|A(t,x)\|^2\le K_1[1+\|x\|^2];$$
$$\|a(t,x_1)-a(t,x_2\|^2+\|A(t,x_1)-A(t,x_2)\|^2\le L_1\|x_1-x_2\|^2;$$
 $$\|c(t,x)h\|^2+\|C(t,x)h\|^2\le K_2\|h\|^2; $$
$$\|c(t,x_1)-c(t,x_2)h\|^2+\|[C(t,x_1)-C(t,x_2)]h\|^2\le L_2\|x_1-x_2\|^2\|h\|^2.$$
Recall that we use notation  $\|A\|=[\sum_{j,k=1}^d A_{kj}A_{jk}]^{\frac 12}$ for $A\in L(R^d)$.

{\bf Lemma 2.1.}
 Let condition  {\bf C 2.1} hold. Then there exists a unique solution
 $\kappa(t)=(\xi(t),\eta(t))\in R^d\times R^{d_1}$ to (\ref{2.1}),(\ref{2.2}) such that $\xi(t)\in R^d$ is a Markov process with  $E\|\xi(t)\|^2<\infty$ and $\eta(t)\in R^{d_1}$  with  $E\|\eta(t)\|^2<\infty$ for any $t\in [0,T]$.

It follows from  {\bf C 2.1} that coefficients of equations  (\ref{2.1}) and (\ref{2.2})  satisfy classical conditions of the existence and uniqueness theorem for solutions of SDEs and hence the lemma statement results from this theorem.

{\bf Lemma 2.2.}  Let condition  {\bf C 2.1} hold. Then the stochastic process $\eta(t)$ satisfying (\ref{2.2}) gives rise to a multiplicative operator functional $\Gamma(t)\equiv \Gamma(t,s):{\cal H}^2_{s}(R^{d_1})\to {\cal H}^2_{t}(R^{d_1})$ satisfying the SDE
\begin{equation}\label{2.3}
d\Gamma(t)= c(t,\xi(t))\Gamma(t)dt+ C(t,\xi(t))(\Gamma(t), dw(t)) ,\quad \Gamma(s,s)=I,
\end{equation}
where $I$ is the identity operator in $R^{d_1}$. Moreover there
exists an inverse map  $\Gamma^{-1}(s,t):{\cal H}^2_{t}(R^{d_1})\to
{\cal H}^2_{s}(R^{d_1})$ satisfying
\begin{equation}\label{2.4}
\Gamma^{-1}(s,t)=I-\int_{s}^{t}\Gamma^{-1}(\theta,t)[c(\theta,\xi(\theta))-C^2(\theta,\xi(\theta))]d\theta-
\int_{s}^{t}\Gamma^{-1}({t})C(\theta,\xi(\theta))dw(\theta)
\end{equation}
with probability 1.

Proof.
%\smartqed
 Under the condition  {\bf C 2.1} we can state the existence and uniqueness of  a solution to  (\ref{2.4}) and the corresponding properties of the map $\Gamma^{-1}(s,t)$. In particular we deduce from uniqueness of  solutions to  (\ref{2.2}) and (\ref{2.4})  that  the  map $\Gamma(t,s)$ defined by
$\eta(t)=\Gamma(t,s)h$ is an evolution family, that is
$\Gamma(t,\theta)\Gamma(\theta,s)=\Gamma(t,s)$  with probability 1
and the map $\Gamma^{-1}(t,s)$ has the same property. Besides by
Ito's formula we can check that $\Gamma(t,s)\Gamma^{-1}(s,t)=I$ a.s.

 Let  $\Gamma^*(s,t)$ be defined by $\langle \Gamma(t,s) h, u\rangle
=\langle h, \Gamma^*(s,t)u\rangle.$ We can verify that
$\Gamma^*(s,t)$ is an invertible evolution map acting from ${\cal
H}^2_{t}(R^{d_1})$ to ${\cal H}^2_{s}(R^{d_1})$.
 Here and below we identify the space $R^d$ with its dual
space $(R^d)^*$.

Consider  a BSDE of the form
\begin{equation}\label{2.5}
dy(t)= -\Gamma^*(s,t)g([\Gamma^*]^{-1}(s,t)y(t),
[\Gamma^*]^{-1}(s,t)z(t))dt+ z(t) dw(t),\, y(T)=\zeta,
\end{equation}
and state conditions on its parameters $g$ and $\zeta$ to ensure
that there exists
  a unique solution $(y(t)\in R^d_1, z(t)\in R^{d\times d_1})$  to  (\ref{2.5}).

We say that condition
 {\bf C 2.2} holds when:

$g:[s,T]\times R^d\times R^{d_1}\times R^{d\times d_1}\to R^{d_1}$,
$\zeta\in R^{d_1}$ be an ${\cal F}_T$-measurable square integrable
random variable and   there exist constants  $L, L_3,$ such that

$$\|g(t,x^1,y,z)-g(t,x^2,y,z)\|\le L _3\|x^1-x^2\|,$$
$$\|g(t,x,y^1,z^1)-g(t,x,y^2,z^2)\|\le L[\,\|y^1-y^2\|+\|z^1-z^2\|\,],$$
$$\langle y-y_1, g(t,x,y^1,z)-g(t,x,y^2,z)\rangle \le \mu\|y-y^1\|^2,$$

3) There exists a constant  $C_0>0$    such that for all $x,x'\in R^d$
$$\|u_0(x)-u_0(x')\|\le C_0\|x-x'\|.$$

 Denote by $f(t,y,z)=\Gamma^*(t)g( \xi(t), [\Gamma^*]^{-1}(t)y, [\Gamma^*]^{-1}(t)z)$ and let $\zeta=\\ \Gamma^*(s,T)u_0(\xi(T))$, where $\xi(t), t\in[s,T]$  is a solution to (\ref{2.1}).  Consider a BSDE
 \begin{equation}\label{2.6}
dy(t)= -f(t,\xi(t), y(t), z(t))dt+ z(t) dw(t) ,\quad y(T)=\zeta\in
R^{d_1}.
\end{equation}

A couple  of progressively measurable random processes $(y(t),z(t))\in {\cal B}^2$
is called a solution of (\ref{2.6}) if with probability 1
\begin{equation}\label{2.7}
y(t)=\zeta+\int_t^Tf(\theta,\xi(\theta),
y(\theta),z(\theta))ds-\int_t^T z(\theta) dw(\theta),\quad 0\le t\le
T.
\end{equation}

{\bf Lemma 2.3.} Let conditions {\bf C 2.1}, {\bf C 2.2} hold. Then
$$\|f(t,x,y_1,z_1)-f(t,x,y_2,z_2)\|\le L[\,\|y_1-y_2\|+\|z-z_1\|\,].$$

Proof.
%\smartqed
 By Lipschitz continuity of $g$  and the properties of $\Gamma(t)$ we    have a.s.
$$\|f(t,y_1,z_1)-f(t,y_2,z_2)\|=\|g(t,\xi(t),[\Gamma^*]^{-1}y_1 , [\Gamma^*]^{-1}z_1)- g(t,\xi(t),[\Gamma^*]^{-1}y_2 , [\Gamma^*]^{-1}z_2)\|$$$$\le \|\Gamma^*\|L[\,\|[\Gamma^*]^{-1}y_1 -[\Gamma^*]^{-1}y_2\|\,]+[\,\|[\Gamma^*]^{-1}z_1-[\Gamma^*]^{-1}z_2\|\,\,]\le L[\,\|y-y_1\|+\|z_1-z_2\|\,].$$

 Given $(u,v)\in{\cal B}^2$,  we define a map  $M$  by $(y,z)=M(u,v)$ as follows.
Let $\zeta$ be $R^{d_1}$-valued ${\cal F}_T$- measurable random
variable and given $f:[s,T]\times R^d\times R^{d_1}\times R^{d\times
d_1}\to R^{d_1}$ set
 \begin{equation}\label{2.8}
y(t)=E[\zeta +\int_t^Tf(\theta,\xi(\theta),
u(\theta),v(\theta))d\theta|{\cal F}_t],\quad 0\le t\le T.
\end{equation}
We apply the Ito theorem about martingale representation of a square
integrable random variable
$$\chi=\zeta+\int_0^Tf(\theta,u(\theta),v(\theta))d\theta$$
to define the process $z(t)$ by the equality
$$\chi=E[\chi]+\int_0^Tz(\theta)dw(\theta).$$
It is easy to check that the couple $(y,z)$ defined in this way satisfies
 $$y(t)= \zeta +\int_t^Tf(\theta,\xi(\theta),u(\theta),v(\theta))d\theta -\int_t^T z(\theta)dw(\theta).$$
In a standard way we show that $M$ acts in ${\cal B}^2$ and  possesses a contraction property.
To this end we denote by
$\bar f=f_1-f_2$ for $f= y,z,u,v$. By Ito's formula we obtain
$$e^{\beta t}E\|\bar y(t)\|^2 +E\left[\int_t^Te^{\beta s}[\beta\|\bar y(s)\|^2+\|\bar z(s)\|^2]ds\right] =
$$$$= 2E\left[\int_t^Te^{\beta s}\langle \bar y(s), f(s, u_1(s), v_1(s))-f(s, u_2(s),v_2(s))\rangle ds\right].$$
Taking into account  Lipschitz continuity of $f$ we obtain
$$E[e^{\beta t}\|\bar y(t)\|^2] +E\left[\int_t^Te^{\beta s}[\beta\|\bar y(s)\|^2+\|\bar z(s)\|^2]ds\right]
\le $$
$$\le 2LE\left[\int_t^Te^{\beta s}\|\bar y(s)\|[ \| \bar
u(s)\|+\|\bar v(s)\|] ds\right]$$ and by the elementary inequality
$2ab\le a^2\alpha^2+\frac{b^2}{\alpha^2},$
$$E[e^{\beta t}\|\bar y(t)\|^2]+E\left[\int_t^Te^{\beta s}\|\bar z(s)\|^2]ds\right] \le
$$$$\le [2L^2 \alpha^2-\beta]E\left[\int_t^T e^{\beta s}\|\bar y(s)\|^2ds\right]+\frac 1{\alpha^2}E\left[
\int_t^Te^{\beta s}(\|\bar u(s)\|^2+\|\bar v(s)\|^2)ds\right]. $$

Choosing $\frac 1{\alpha^2}=\frac 12$ and $\beta -4L^2= 1$ we obtain
$$e^{\beta t}E\|\bar y(t)\|^2 +E\left[\int_t^Te^{\beta s}\|\bar z(s)\|^2ds\right] \le
\frac 12 E\left[\int_t^Te^{\beta s}[\|\bar u(s)\|^2+\|\bar v(s)\|^2]ds\right]. $$
In the similar way we can check that  $(y,z)=M(u,v)\in{\cal B}^2$.
As a result we deduce that $M$ is a contraction in ${\cal B}^2$ and the following statement holds.

{\bf Theorem 2.1}   Let   condition {\bf C 2.2} hold. Then there exists a unique solution  $(y,z)\in{\cal B}^2$  of BSDE  (\ref{2.6}) and successive approximations $(y^n,z^n)$ of the form
$$y^{n+1}(t)=\zeta+\int_t^Tf(\theta, \xi(\theta),y^n(\theta), z^n(\theta))d\theta-\int_t^Tz^{n+1}(\theta)dw(\theta)$$
converges to  the solution of (\ref{2.6}) with probability 1.

Proof.
%\smartqed
  The existence and uniqueness of a solution $(y,z)$ to  (\ref{2.6}) follows from the fixed point theorem for the  contraction  $M:{\cal B}^2\to{\cal B}^2$.  Applying the above estimates to the successive approximations $(y^n,z^n)$
we can verify that
$$E\left[\int_t^Te^{\beta\theta}\| y^n(\theta)-y^m(\theta)\|^2 ds|{\cal F}_t\right] +
E\left[\int_t^Te^{\beta\theta}\| z^n(\theta)-z^m(\theta)\|^2
ds|{\cal F}_t\right] \to 0,  \quad m,n\to\infty$$ with probability
1. Hence, $(y^n,z^n)$ is a Cauchy sequence in ${\cal B}^2$
 and the limit $P-\lim_{n\to\infty}(y^n,z^n)=(y,z)$ exists and satisfies (\ref{2.4}).

 Below along with  a weakly coupled multidimensional FBSDE  of the form
 \begin{equation}\label{2.9}
dy(t)= -f(t,\xi(t),y(t),z(t))dt+z(t)dw(t),\quad y(T)=\Gamma^*(s,T)u_0(\xi(T)),
\end{equation}
where $\xi(t)$ is a solution of  (\ref{2.4}) we consider a
 weakly coupled scalar  FBSDE   which can be described as follows.
 Let
 \begin{equation}\label{2.10}q(\kappa)=\pmatrix{a (x)\cr c(x)h},
  Q(\kappa)=\pmatrix{A(x)&0\cr0&C(x)h},\quad \tilde G(\kappa,y,z)=\langle h, f(x,y,z)\rangle.\end{equation}

Obviously,   we can rewrite the system (\ref{2.1}),(\ref{2.2}) in  the form
 \begin{equation}\label{2.11}
d\kappa(t)= q(t,\kappa(t))dt+Q(t,\kappa(t))dW(t),\quad \kappa( s)= \kappa=(x,h),
\end{equation}
The required  FBSDE  can be presented in  the form
 \begin{equation}\label{2.12}
dY(t)= -\tilde G(t,\kappa(t),Y(t),Z(t))dt+\langle
Z(t),dW(t)\rangle,\quad Y(T)=\langle \eta(T), u_0(\xi(T)\rangle,
\end{equation}
where $\kappa(t)=(\xi(t),\eta(t))$ solves (\ref{2.11}),
$W(t)=(w(t),w(t))^*$ and $\langle Z(t),dW(t)\rangle= \langle h,
z(t)dw(t)\rangle.$

A triple of progressively measurable random processes $(\kappa(t),y(t),z(t))\in {\cal B}^3 $
is called a solution of (\ref{2.11}),(\ref{2.12}) if with probability 1 for  all  $0\le s\le t\le T$
\begin{equation}\label{2.13}
\kappa(t)=\kappa+\int_s^t q(\theta,\kappa(\theta))d\theta+
\int_s^tQ(\theta,\kappa(\theta))dW(\theta),
\end{equation}
\begin{equation}\label{2.14}
Y(t)=\langle \eta(T), u_0(\xi(T))\rangle +\int_t^T\tilde
G(\theta,\kappa(\theta),Y(\theta),Z(\theta))d\theta-\int_t^T \langle
Z(\theta), dW(\theta)\rangle.
\end{equation}
The FBSDEs (\ref{2.1}), (\ref{2.2}),  (\ref{2.6})   and
(\ref{2.11}), (\ref{2.12}) are equivalent.

\section{Comparison theorem for multidimensional BSDE}
\label{sec:3}
\setcounter{equation}{0}
Comparison theorems present an important tool in the BSDE and FBSDE theory and in particular
in the context of the connections between  FBSDE theory and
viscosity solutions of corresponding parabolic equations and systems.  In this paper to prove a
comparison theorem for a multidimensional   BSDE  we use the special features of the BSDE under consideration.

Consider a couple of $d_1$-dimensional BSDEs
\begin{equation}\label{3.1}
y^i(t)=\zeta^i+\int_t^Tf^i(\theta,y^i(\theta),z^i(\theta))d\theta-\int_t^T
z^i(\theta) dw(\theta),\quad i=1,2
\end{equation}
for $0\le t\le T$  and  use the specific  features of these BSDEs
investigated in the previous sections. Here
$\zeta^i,f^i(\theta,y,z)\in R^{d_1} $ for $\theta\in[0,T], y\in
R^{d_1}, z\in R^{d\times d_1}.$

For any fixed nonzero vector $h\in R^{d_1}$ and $y^1, y^2\in
R^{d_1}$ we say that $y^1\le_h y^2$ under $h$ if $\langle
h,y^1\rangle\le \langle h,y^2\rangle$. Without loss of generality we
choose $h$ to have $\|h\|=1$.

Given two vectors $y^1,y^2\in R^{d_1}$, we say $y^1\le y^2$ if
$y_m^1\le y^2_m$,$m=1,\dots , d_1$, where $y_m=\langle y,e_m\rangle$
and $(e_m)_{m=1}^{d_1}$ is a fixed orthonormal basis in $R^{d_1}$.

Given $f\in R^{d_1}$ we denote by  $f^+_m=\max[f_m,0], \, m=1,\dots, d_1.$

Consider a couple of BSDEs with parameters $\zeta^i,  f^i$, $i=1,2$.

We say that condition {\bf C 3.1} holds if

i) $\zeta^1\le \zeta^2, \, P-\mbox{a.s.}\quad, $

ii) for each $m=1,\dots, d_1$ inequality  $f^1_m(t,y^1,z^1)\le
f^2_m(t,y^2, z^2)$  holds true  when $y^1_l\le y^2_l$ for all
$l=1,\dots, d_1$ except $l=m$    while  $y^1_m=y^2_m,$ and $ \,
z_{mk}^1=z_{mk}^2$ for each   $k=1,\dots, d$,

iii) For all $y^1,y^2\in R^{d_1}, z^1, z^2\in R^{d\times d_1}$ and
for each $m=1,\dots,d_{1}$
$$\|f_m^i(t,y^1,z^1)-f^i_m(t,y^2,z^2)\|\le L[\|y^1-y^2\|+\|z^1-z^2\|], i=1,2.$$

Set   $\bar \alpha=\alpha^1-\alpha^2$   for $\alpha= y,\zeta, f$ and $ z$ as well.

Let us mention that within this section we do not assume summing up with respect to repeating indices.

{\bf Theorem 3.1.}     Let  $(\zeta^i, f^i),\, i=1,2$ be parameters of BSDEs  (\ref{3.1}) satisfying
conditions  {\bf C 2.1} and   {\bf C 3.1}.   Assume  that   $(
y^i(t),z^i(t)),\, i=1,2,\quad  t\in[s,T]$ solve   (\ref{3.1})  with
this parameters. Then  $y^1(t)\le y^2(t)$ a.s.  Moreover the
comparison is strict, that is if in addition $y^2(s)=y^{1}(s)$ then
$\zeta^1=\zeta^2$, $f^2(t,y^2(t),z^2(t))=f^{1}(t,y^{2}(t),z^{2}(t))$
and $y^2(t)=y^{1}(t),\,\forall t\in[s,T] P$-a.s. In particular
whenever either  $P(\zeta^1<\zeta^{2})>0$ or $
f^1(t,y^2(t),z^2(t))<f^{2} (t,y^2(t),z^2(t))$ on a set of positive
$dt\times dP$ measure, then $y^1(s)< y^{2}(s)$ a.s.

Proof.
%\smartqed Set $\bar y_j(t)^+=\max[y_j^1(t)-y_j^2(t),0]$.
Applying Ito's formula to $ |\bar y_j(t)^+|^2$ where $j=1,\dots,
d_1$, and evaluating mean value
 we get
\begin{equation}\label{3.2}
E|\bar y_j(t)^+|^2= E|\bar\zeta_j^+|^2-E[\int_t^T2I_{y_j^1(s)>y^2_j(s)}\bar y_j(s)[f_j(s,y^1(s),z^1(s))-\end{equation}$$f_j(s,y^2(s),z^2(s))]ds]
-E\left[\int_t^T I_{\{ y_j(s)>y_j^2(s)\}}\|\bar z_j(s)\|^2ds\right] -E\left[\int_t^T \bar y_j^+dL_j(s)\right],
$$
where $L_j(t)$ is the local time of $\bar y_j(s)$  at 0.  Note that
the last summand is equal to 0 and since
 $\zeta^1\le \zeta^2$ a.s. we have $E[\|[\zeta^1-\zeta^2]^+\|^2]=0$,    Obviously,\\
$E\left[\int_t^TI_{y_j^1(s)>y^2_j(s)}\bar y_j(s)\bar z_j(s)dw(s)\right]=0.$
Hence,
$$E[\bar y_j(t)^+]= E\left[\int_t^TI_{y_j^1(s)>y^2_j(s)}2 \bar y^+_j(s)[f_j^1(s,y^1(s), z^1(s))-
f_j^2(s,y^2(s), z^2(s))]ds\right]-$$$$E\left[\int_t^T I_{\{
y^1_j(s)>y_j^2(s)\}}\|\bar z_j(s)\|^2ds\right].$$

Set
$$\bar f_j(s)=f_j^1(s,y^1, z^1)-f_j^2(s,y^2, z^2)=
$$$$f_j^1(s,y^1_1,\dots, y_j^1,\dots,y_{d_1}^1, z^1_1,\dots, z_j^1,\dots, z^1_{d_1})-$$
$$-f_j^2(s,y^2_1,\dots, y^2_j, \dots, z^2_{1},\dots, z_j^2,\dots,z_{d_1}^2)=$$
$$
=[f_j^1(s,y^1_1,\dots, y_j^1,\dots,y_{d_1}^1, z^1_1,\dots,
z_j^1,\dots, z^1_{d_1})-$$$$ f_j^2(s,y^1_1+\bar y_1^+,\dots,
y_j^1,\dots,y_{d_1}^1+\bar y^+_{d_1}, z_1^2,\dots, z_j^1,\dots,
z^2_{d_1})]+
$$
$$
+[f_j^2(s,y^2_1+\bar y_1^+,\dots, y_j^1,\dots,y_{d_1}^2+\bar
y^+_{d_1}, z^2_1,\dots, z_j^1,\dots, z^2_{d_1})-$$
$$
-f_j^2(s,y^2_1,\dots, y^2_j, \dots, z^2_{1},\dots,
z_j^2,\dots,z_{d_1}^2)]=\Pi_1+\Pi_2
$$

Since for any  $m =1,\dots,d_1$ we have $y^1_m\le y^2_m+\bar y^+_m$
for $m\ne j$, taking into account  ii) in {\bf C 3.1}  we get
$\Pi_1\le 0$.

Next, due to Lipschitz continuity of $f^2$ we  have
$$\Pi_2\le L[|\bar y^+_1|+\dots+|\bar y_{j-1}^+|+|\bar y_j|+\dots+|\bar y^+_{d_1}|+\|\bar z_j\|].$$

Applying Ito's formula  due to generator properties  we deduce that
$$E|\bar y^+_j(t)|^2\le 2E\left[\int_t^TI_{y_j^1(s)>y^2_j(s)}\bar
y_j^+(s)\bar f_j(s)ds\right] -
E\left[\int_t^TI_{y_j^1(s)>y^2_j(s)}\sum_{k=1}^d|\bar
z_{jk}(s)|^2ds\right]\le$$
$$\le E\left[2\int_t^TI_{y_j^1(s)>y^2_j(s)}L\bar y_j^+(s)[|\bar
y_1(s)|+\dots+ |\bar y_{j-1}^+|+ |\bar y_j(s)|+\dots+|\bar
y^+_{d_1}|+\|\bar z_j(s)\|]ds\right]-
$$
$$ -E\left[\int_t^T I_{\{ y_j^1(s)>y_j^2(s)\}}\|\bar z_j(s)\|^2ds\right]\le
 E\left[\int_t^T I_{\{ y_j^1(s)>y_j^2(s)\}} L^2(d_1+1)|\bar y_j(s)|^2
 ds\right]+$$
 $$+ E\left[\int_t^T I_{\{ y_j^1(s)>y_j^2(s)\}}[\sum_{k=1}^{d_1}|\bar y_k(s)|^2 +\|\bar z_j(s)\|^2]ds\right]-$$$$E\left[\int_t^T I_{\{ y_j^1(s)>y_j^2(s)\}}\|\bar
z_j(s)\|^2ds\right]= L^2(d_1+1)\int_t^TE [I_{\{
y_j^1(s)>y_j^2(s)\}}|\bar y_j(s)|^2 ]ds$$
\begin{equation}\label{3.04}+\int_t^TE [I_{\{ y_j^1(s)>y_j^2(s)\}}\sum_{k=1}^{d_1}|\bar y_k(s)|^2 ]ds.\end{equation}
Note that above we have used  an elementary inequality of the form
$$2L\bar y^+_j(s)|\bar y_k(s)|\le L^2|\bar y_j^+(s)|^2+|\bar y_k(s)|^2.$$

 Summing up left and right hand side in (\ref{3.04})  we get that the function
 $m(t)=\sum_{j=1}^{d_1}E|\bar y^+_j(t)|^2$
 satisfies inequality
$$m(t)\le (L^2(d_1+1)+d_1)\int_t^Tm(s)ds$$
    Finally,  due to results of the previous section we know that for
 $t\in[0,T]$ the inequality  $E|\bar y_j(t)^+|^2<\infty$  holds for each $j=1,\dots, m$ then
 by the Gronwall lemma we know that $m(t)=0$ and since
    $m$ is a sum of positive summands, each summand should be equal
    to zero.
  Hence $|\bar y_j^+(t)|= 0$ and thus $y^1_j(t)\le
y^2_j(t)$ a.s. for all $j=1,\dots, d_1.$

At the end of this section we come back to the one-dimensional BSDE  (\ref{2.14})  and derive the corresponding  comparison theorem. Note that this theorem motivates  our choice of comparison for vector functions in the case under consideration.

Consider the SDE of the form
\begin{equation}\label{3.4}
\kappa(t)=\kappa+\int_s^t
q(\kappa(\theta))d\theta+\int_s^tQ(\kappa(\theta))dW(\theta), \quad
 s\le t\le T,
\end{equation}
introduced in the previous section and note that  one can consider instead of the  BSDE
\begin{equation}\label{3.5}
y(t)=\Gamma^*(s,T)u_0(\xi(T))+\int_t^Tf(\theta,\xi(\theta),y(\theta),z(\theta))d\theta-\int_t^T
z(\theta) dw(\theta),\quad s\le t\le T,
\end{equation}
  with respect to the process $y(t)\in R^{d_1}$     a new BSDE
 \begin{equation}\label{3.6}
dY(t)=- \tilde G(t,\kappa(t),Y(t),Z(t))dt+ \langle
Z(t),dW(t)(t)\rangle,\quad Y(T)=\zeta= \langle \eta(T),u_0(\xi(T))\rangle,
\end{equation}
where  $Y(t)=\langle \eta(t), u(t,\xi(t))\rangle$ is a scalar
process. We denote  $|Y|=\sup_{\|h\|=1}| \langle h,u\rangle|=
\|u\|$.

{\bf Theorem 3.2.}
Let $(Y^i,Z^i), i=1,2$ be solutions of one dimensional BSDEs

\begin{equation}\label{3.7}
dY^i(t)=- \tilde G^i(t,\kappa(t),Y^i(t),Z^i(t))dt+ \langle
Z^i(t),dW(t)\rangle,\quad Y^i(T)=\Upsilon^i= \langle
\eta(T),u^i_0(\xi(T))\rangle.
\end{equation}

Suppose that $\Upsilon^1\le \Upsilon^2$ and $\tilde G^1(t,\kappa,
Y^2,Z^2)\le \tilde G^2(t,\kappa, Y^2,Z^2)$ $dt\times dP$ - a.e. Then
$Y^1(t)\le Y^2(t)$  a.s. for all $s\le t\le T$.

Proof
%\smartqed
 Define a scalar process
$$
\mu(t)=\cases{\frac{\tilde G^1(t,\kappa(t),Y^2(t),Z^1(t))-\tilde
G^1(t,\kappa(t),Y^1(t),Z^1(t))}{Y^2(t)-Y^1(t)}& if $Y^1(t)\neq
Y^2(t),$\cr 0& if $Y^1(t)= Y^2(t),$}$$ and a vector process
$\nu(t)\in R^{d}$ such that
$$\nu_k(t)=\cases{\frac{\tilde
G^1(t,\kappa(t),Y^1(t),Z^{(k)}(t))-\tilde
G^1(t,\kappa(t),Y^1(t),Z^{(k-1)}(t))}{Z^2_k(t)-Z^1_k(t)}& if
$Z^1_k(t)\neq Z^2_k(t)$\cr 0& if $Z^1_k(t)= Z^2_k(t)$},$$ where
$Z^{(k)}(t)$ denotes the $d$-dimensional vector such that its first
$k$ components  are equal to corresponding components of $Z^2$   and
the remaining $d-k$ components are equal to those of $Z^1$. Due to
Lipschitz continuity of $g$ the processes $\mu(t)$ and $\nu(t)$ are
bounded  and in addition they are progressively measurable.

 As above we use notation $\bar f=f^1-f^2$ for $f=Y,Z,\Upsilon$ and observe that $(\bar Y(t),\bar Z(t)) $ satisfies the BSDE
 $$\bar Y(t)=\bar \Upsilon+\int_t^T[\mu(\theta)\bar Y(\theta)+
 \langle \nu(\theta),\bar Z(\theta)\rangle]d\theta+\int_t^TN(\theta)d\theta-\int_t^T\langle \bar Z(\theta),
  dW(\theta)\rangle,$$
 where $N(t)=\tilde
G^1(t,\kappa(t),Y^2(t),Z^2(t))-\tilde
G^2(t,\kappa(t),Y^2(t),Z^2(t)).$ For $s\le t\le T$ we define
$$\rho_{s,t}=exp\left[\int_s^t(\mu(\theta)-\frac 12\|\nu(\theta)\|^2)d\theta +
\int_s^t\langle\nu(\theta), dW(\theta)\rangle\right].$$ By Ito's
formula we can verify that $(\bar Y(\theta),\bar Z(\theta))$ satisfy
the BSDE
$$ d[\rho_{s,\theta}\bar Y(\theta)]=\rho_{s,\theta}[\bar Y(\theta)+N(\theta)]d\theta+
\rho_{s,\theta}\langle \bar Z(\theta)+ \bar
Y(\theta)\nu(\theta),dW(\theta)\rangle$$ for $\theta\in[s, T]$ and
$$\bar Y(\theta)=E\left[\rho_{s,T}\bar \Upsilon+\int_\theta^T\rho_{s,\vartheta}N(\vartheta)d\vartheta|{\cal F}_\theta\right]$$

The required assertion immediately follows from negativity of $\bar
\Upsilon$ and $N(t)$.

Let us mention a useful remark.   Let  $Y^1, Z^1$ be a solution of
BSDE
$$Y^1(t)=\Upsilon^1+\int_t^T\tilde
G^1(\theta,Y^1(\theta),Z^1(\theta))d\theta-\int_t^T \langle
Z^1(\theta),dW(\theta)\rangle$$ and $(Y^2,Z^2)$ satisfy
$$Y^2(t)=\Upsilon^2+\int_t^TM(\theta)d\theta-\int_t^T \langle Z^2(\theta),dW(\theta)\rangle,$$
where $M(\theta)$ is a scalar progressively ${\cal
F}_\theta$-measurable process. Suppose that $\Upsilon^1\le
\Upsilon^2$ and $\tilde G^1(t,Y^2(t),Z^2(t))\le M(t)$.  Then we can
choose
$$\tilde
G^2(t,\kappa(t),Y^{2},Z^{2})=\tilde
G^1(t,\kappa(t),Y^{2},Z^{2})+[M(t)-G^1(t,\kappa(t),Y^2(t),Z^2(t))]$$
and apply the result of theorem 3 to deduce that $Y^1(t)\le Y^2(t)$.
If in addition $\tilde G^1(t,\kappa(t),Y^2,Z^2)<M(t)$ on a set of
positive measure $dt\times dP$, then $Y^1(s)< Y^2(s)$.

\section{Viscosity solution to nonlinear parabolic system}
\label{sec:4}

\setcounter{equation}{0}
In this section we show that a solution of a forward-backward stochastic differential equation
generates a viscosity solution of the Cauchy problem for  a system of quasilinear parabolic equations.

 Let $(\xi(t)\in R^d, y(t)\in R^{d_1},z(t)\in R^{d\times d_1})$ be a solution of the FBSDE

 \begin{equation}\label{4.1}
d\xi(t)=a(\xi(t))dt +A(\xi(t))dw(t),\quad \xi(s)=x,
\end{equation}
\begin{equation}\label{4.2}
dy(t)= -\Gamma^*(t)g([\Gamma^*]^{-1}(t)y(t),
[\Gamma^*]^{-1}(t)z(t))dt+ z(t) dw(t) ,\quad 
\end{equation}
$$y(T)=\Gamma^*(s,T)u_0(\xi(T)),$$
where $\Gamma(t)$ is a multiplicative operator functional of the
process $\xi(t)$ generated by the solution $\eta(t)\in R^{d_1}$ of
the linear SDE
 \begin{equation}\label{4.3}
d\eta(t)=c(\xi(t))\eta(t)dt +C(\xi(t))(\eta(t),dw(t)),\quad
\eta(s)=h,
\end{equation}
and $u_0:R^d\to R^{d_1}$ be a continuous bounded function.

Denote by $S^{d_1}_+=\{h\in R^{d_1}: h_m\ge
0,m=1,\dots,d_1\,\mbox{and }\,\|h\|=1\},$ and let $e_1,\dots,
e_{d_1}$ be a fixed orthonormal basis in $R^{d_1}$.

In section 2 we have shown that one can write (\ref{4.2}) in the form
\begin{equation}\label{4.4}
dy(t) = -f(t,\xi(t),y(t), z(t))dt+ z(t) dw(t), \quad
y(T)=\Gamma^*(s,T)u_0(\xi(T)),
\end{equation}
and proved that given a solution $\xi(t)$  of  (\ref{4.1}),
there exists a unique solution $(y(t),z(t))$ of this BSDE.

  Assume that there exists a solution $(\xi_{s,x}(t), y^{s,x}(t),z^{s,x}(t))$   to  (\ref{4.1}), (\ref{4.2})
    and the comparison theorem 2  is valid.  The aim of this section is to prove that the function  $u(s,x)=y^{s,x}(s)$ is a viscosity solution of the Cauchy problem
  \begin{equation}\label{4.5}
\frac{\partial u_l}{\partial s}+\frac 12 Tr A^*(x)\nabla^2
u_lA(x)+\langle a(x),\nabla u_l\rangle+
\end{equation}
$$
+B_{lm}^i(x)\nabla_i u_m+c_{lm}(x)u_m +g_l(x,u, K(u,\nabla
u))=0,\quad l=1,\dots,d_1,
$$
$$ u(T,x)=u_0(x),$$
where $B_{lm}^i=\sum_{q=1}^dC^q_{lm}A^{qi},\quad K(u,\nabla
u)=C^*u+A^*\nabla u.$

As it was mentioned in section 2 the system (\ref{4.5})   can be
easily reduced to a scalar parabolic equation
 \begin{equation}\label{4.6}
\frac{\partial V}{\partial s}+\frac 12 Tr Q^*(x,h)\nabla^2 V
Q(x,h)+\langle q(x,h),\nabla V\rangle +G(h, x,V, Q^*\nabla
V)=0,\end{equation}
$$ V(T,x)=V_0(x,h)=\langle h,u_{0}(x)\rangle $$
with respect  to a scalar function $V$ defined on $[0,T]\times
R^d\times S^{d_1}_+$  (see equation (\ref{2.1})).

Hence we recall first the  definition of a viscosity solution of the
Cauchy problem for a general scalar nonlinear
 parabolic equation
\begin{equation}\label{4.7}\frac{\partial V}{\partial s}+\Psi(s,z,V,\nabla V,\nabla^2V)=0.\quad V(T,z)=V_0(z),\end{equation}
where $z=(x,h).$

 A function  $\Psi: [0,T]\times (R^d\times S^{d_1}_+)\times R\times
(R^d\times R^{d_1})\times R^d\otimes R^d\to R$ satisfying estimates
$$\Psi(s,z,V,p,q)\le \Psi(s,z,U,p,q) \quad \mbox{if} \quad V\le U,$$
and
$$\Psi(s,z,V,p,q)\le \Psi(s,z,V,p,q_1) \quad \mbox{if} \quad q_1\le q$$
is called  a proper function.

Given a proper function $\Psi$  to define a viscosity solution of
(\ref{4.7})  one has to introduce notions of a sub- and  a
supersolution of this Cauchy problem.

Denote by $C^{1,2}_{d,d_1}\equiv C^{1,2}([0,T]\times R^d; R^{d_1})$
a set of functions $\psi:[0,T]\times R^d; R^{d_1}$ differentiable in
$s\in[0,T]$ and twice differentiable in $x\in R^d$.

A continuous real valued function $V(s,z)$ is called  a subsolution
of (\ref{4.7}) if $V(T,z)\le V_0(z)$, $z\in R^{d_2} ,\, d_2=d+d_1, $
and for any $\Phi\in C^{1,2}_{d_2,1}$ and a point $(s,z)\in
[0,T]\times R^{d_2}$ which is a local maximum of $V(t,\tilde
z)-\Phi(t,\tilde z)$ the inequality
$$\frac{\partial \Phi}{\partial s}+\Psi(s,z,V,\nabla \Phi,\nabla^2\Phi)\ge 0$$
holds.

 A continuous function $V(s,z)$ is called  a super-solution of
(\ref{4.7}) if $V(T,z)\ge V_0(z)$, $z\in R^{d_2}$ and for any
$\phi\in C^{1,2}_{d_2,1}$ and $(s,x)\in [0,T]\times R^d$   which is
a local minimum of $u_m(t,\tilde x)-\phi_m(t,\tilde x)$ the
inequality
$$\frac{\partial \Phi}{\partial s}+\Psi(s,z,V,\nabla \Phi,\nabla^2\Phi)\le 0$$
holds. A continuous function $V(s,z)$ is called  a viscosity
solution of (\ref{4.7}), if it is both sub- and super-solution of
this Cauchy problem. Hence to prove that the function  $V(s,z)$  is
a viscosity solution to  (\ref{4.7}) one has to prove that  $V$
 is both  sub- and  supersolution of  (\ref{4.7}).

To give a definition of a viscosity solution of the Cauchy problem
to the system (\ref{4.5}) we use a definition of a viscosity
solution of the scalar Cauchy problem (\ref{4.6}) and then rewrite
the definition in terms of the solution to (\ref{4.5}).

Given functions  $\phi_m\in C^{1,2}_{d,d_1}$, $m=1,\dots,d_1$ denote
by
$$[{\cal A}\phi]_m(x)=\frac 12 Tr A^*(x)\nabla^2 \phi_mA(x)+\langle a(x),\nabla \phi_m\rangle
+B_{ml}^i(x)\nabla_i \phi_l+c_{ml}(x)\phi_l,$$ where $i=1,\dots d,\,
m,l=1,\dots, d_1$.

Let $(s,x,\phi,p,q)\in [0,T]\times R^d\times R^{d_1}\times R^{d\times d_1}\times R^{d^2\times d_1}$ and
\begin{equation}\label{4.8}{\cal M}_m(s,x,\phi,p,q_m)=\frac 12 Tr A^*(x)q_mA(x)+\langle
a(x),p_m\rangle + \end{equation} $$ +B_{ml}^i(x)\nabla_i
p_l+c_{ml}(x)\phi_l +g_l(s,x,u, p).
$$
Given ${\cal M}_m, m=1,\dots, d_1,$ of the form (\ref{4.8}) the
system
\begin{equation}\label{4.9}
\frac{\partial u_m}{\partial  s}+{\cal M}_m(s,x,u,\nabla
u,\nabla^2u_m)=0
\end{equation}
coincides with (\ref{4.5}).

A
 continuous  function $u:[0,T]\times R^d \to  R^{d_1}$ is called  a sub-solution of  (\ref{4.9}) if
  for each $m=1,\dots,d_1$ an inequality
 $u_m(T,x)\le u_{0m}(x)\rangle$,   holds
and for any $\varphi_m\in C^{1,2}_{d,1}$ and a point $( s, x)\in
[0,T]\times R^d$ which is a local maximum of
 $u_m(\tilde s,\tilde x)-\varphi_m(\tilde s,\tilde x)$  an inequality
\begin{equation}\label{4.01}\frac{\partial \varphi_m}{\partial  s}+
{\cal M}( s,x,u,\nabla \varphi,\nabla^2\varphi_m)\rangle\ge
0\end{equation} holds.

 A continuous function $u(s,x)$ is called  a super-solution of
(\ref{4.9}) if for each $m=1,\dots,d_1$ an inequality $u_m(T,\tilde
x)\rangle\ge u_{0m}(\tilde x)$, $x\in R^d $ holds and for any
$\varphi_m\in C^{1,2}_{d,1}$ and a point $(s, x)\in [0,T]\times R^d$
which is a local minimum of $u_m(\tilde s,\tilde x)-\varphi_m(\tilde
s,\tilde x)$ an inequality
\begin{equation}\label{4.02}\frac{\partial \varphi_m}{\partial \tilde s}+{\cal M}_m( s,x,u,\nabla
\varphi,\nabla^2\varphi_m)\rangle\le 0,\end{equation} holds.

 A continuous function $u(s,x)$ is called  a viscosity
solution of (\ref{4.9}), if it is both sub- and super-solution of
this Cauchy problem. Hence to prove that the function  $u(s,x)$  is
a viscosity solution to  (\ref{4.9}) one has to prove that  $u$   is
both sub- and super-solution of  (\ref{4.9}).

{\bf Theorem 4.1.} Assume that conditions of theorem 2 hold  and   $(\xi_{s,x}(t), y^{s,x}(t),z^{s,x}(t),\eta^{s,x}(t))$
is a solution to (\ref{4.1})-(\ref{4.3}).  Then
 $u(s,x)=y^{s,x}(s)$  is a continuous in $(s,x)$  viscosity solution of (\ref{4.5}).

Proof.
%\smartqed
  Under assumptions of section 2 continuity of    $u(s,x)=y^{s,x}(s)$  in spatial variable $x$ and time variable $s$ is granted by the BSDE
  theory results [5]  which state that under {\bf C 2.1} and {\bf C 2.2} the solution of BSDE  (\ref{4.4})
  is continuous with respect to   parameters $(s,x)$.  To verify that  $u(s,x)$ is a viscosity solution of
   (\ref{4.5}),  we have to prove that  $u$ is both a subsolution and a supersolution of   (\ref{4.5}).
   First we check that $u$ is  a subsolution.
To this end for each $m=1,\dots, d_1$ we can choose a  function
$\phi_m\in C^{1,2}_{d,1}$ and a point  $(s,x)\in [0,T]\times R^d$
such that
 at the point $(s,x)$ a function $
u_m(s,x)-\phi_m(s,x) $ has a local maximum. Without loss of
generality  we assume that
  $u_m(s,x)=\phi_m(s,x)$.

We have to prove that (\ref{4.01}) holds.

 Assume on the contrary
that there exists $m\in\{1,\dots, d_1\}$ such that
\begin{equation}\label{4.10}
{\cal K}^{u,\phi}_m(s,x)=\frac{\partial \phi_m}{\partial s}+[{\cal
A}\phi]_m(s,x)+ g_m(s,x,u(s,x),K(u, \nabla \phi)(s,x))\rangle < 0.
\end{equation}

By continuity there exists $0<\alpha\le T-s$   such that for all
$\theta\in[s,s+\alpha]$, $x_1\in R^d, h_1\in R^{d_1}$, $\|x-x_1\|\le
\alpha, \|e_m-h_1\|\le \alpha$   the inequalities
\begin{equation}\label{4.03}
\Phi^u(\theta,x_1,h_1)-\Phi^\phi(\theta,x_1,h_1)\le 0\end{equation}
and \begin{equation}\label{4.10} \langle h_1,\left(\frac{\partial
\phi}{\partial \theta}+{\cal A}\phi\right)(\theta,x_1)+
g(\theta,x_1,u(\theta,x_1),K(u, \nabla \phi)(\theta,x_1))\rangle < 0
\end{equation}
hold.

Given $(\xi_{s,x}(t),  \eta_{s,h}(t))$ satisfying (\ref{4.1}),
(\ref{4.3}), define $\tau$ by
 $$\tau=\inf\{t\ge s:\|\xi_{s,x}(t) -x\|\ge
\alpha\}\wedge \inf\{t\ge s:\|\eta_{s,h}(t) -h\|\ge \alpha\}\wedge
(s+\alpha).$$

It follows from results in [10],[11] that the pair
$$(\hat y(t),\hat z(t))= (y^{s,x}(t\wedge \tau), I_{[s,\tau]}(t) z^{s,x}(t\wedge\tau)),\quad s\le t\le s+\alpha$$
satisfies BSDE
\begin{equation}\label{4.11}
\hat
y(t)=\Gamma^{*}(t,\tau)u([s+\alpha]\wedge\tau,\xi([s+\alpha]\wedge\tau))+
\int_t^{s+\alpha}I_{[s,\tau]}(\theta)
f(\theta,\xi(\theta),\hat y(\theta),\hat z(\theta))d\theta-
\end{equation}$$
\int_t^{s+\alpha}   \hat z(\theta)dw(\theta)\rangle,\quad s\le t\le
s+\alpha.$$

On the other hand applying Ito's formula we obtain that the couple
 $$(\tilde y(t),\tilde z(t))=(\Gamma^{*}(t,t\wedge \tau)\phi(t\wedge\tau,\xi_{s,x}(t\wedge \tau)),
  I_{[s,\tau]}(t) K(u,\nabla\phi)(t,\xi_{s,x}(t))),\quad s\le t\le s+\alpha,$$
 where
 $$K(u,\nabla\phi)(t,\xi_{s,x}(t))=\pmatrix{\Gamma^*(t)A^*(\xi(t))\nabla\phi
 (t,\xi_{s,x}(t))\cr
 \Gamma^*(t)C^*(\xi_{s,x}(t))u(t,\xi_{s,x}(t))\cr },\quad s\le t\le s+\alpha,$$
satisfies a  BSDE
$$\tilde y(t)=(\Gamma^{*}( \tau)\phi(\tau,\xi_{s,x}(\tau))+\int_t^{s+\alpha}I_{[s,\tau]}(\theta)
\left(\frac{\partial\phi}{\partial \theta}+[{\cal
A}\phi]\right)(\theta, \xi_{s,x}(\theta))d\theta+
$$$$\int_t^{s+\alpha} \tilde z(\theta)dw(\theta).$$
 Notice that  $ \hat y_m(s)=\tilde
 y_m(s)=u_m(s,x).$

Then for any stopping time $\tau\in [s,s+\alpha]$ due to
(\ref{4.03}) and (\ref{4.10})
 we derive
$$0\ge [ \Phi^u(\tau,\kappa(\tau))-\Phi^\phi(\tau,\kappa(\tau))]=\langle e_m, u(s,x)-\phi(s,x)\rangle-
$$
$$
-\int_s^\tau \langle e_m, [\frac{\partial\phi}{\partial\theta}+{\cal
A}\phi](\theta,\xi_{s,x}(\theta))d\theta-\int_s^\tau \langle e_m,
f(\theta,\xi_{s,x}(\theta),\hat y(\theta),\hat z(\theta))d\theta +
$$
$$
+\int_s^\tau \langle e_m,[\hat z(\theta) - K(u,
\nabla\phi)(\theta,\xi_{s,x}(\theta))] dw(\theta)\rangle.
$$
Keeping in mind that  by assumption for each $m=1,\dots, d_1$ at the
point $(s,x)$  we have $u_m(s,x)-\phi_m(s,x)=0$ and computing the
expectation of both parts of the last inequality we deduce
\begin{equation}\label{4.17}E\left( \int_s^\tau \langle e_m, [\frac{\partial\phi}{\partial\theta}+{\cal
A}\phi](\theta,\xi(\theta))\rangle d\theta+\int_s^\tau \langle
e_m,f(\theta,\xi(\theta),\hat y(\theta),\hat z(\theta))\rangle
d\theta\right)\ge 0.
\end{equation}
Denote by
$$
\gamma_1(s,\tau)=\langle e_m,\int_s^\tau
\{[\frac{\partial\phi}{\partial\theta}+{\cal
A}\phi](\theta,\xi(\theta))+ g(\theta,x_1,u(\theta,\xi(\theta)),K(u,
\nabla \phi)(\theta,\xi(\theta))\rangle,$$
$$\gamma_2(s,\tau)=\langle e_m,\int_s^\tau [f(\theta,\xi(\theta),\tilde y(\theta),
\tilde z(\theta))-g(\theta,\xi(\theta),u(\theta,\xi(\theta)),K(u,
\nabla \phi)(\theta,\xi(\theta))) d\theta\rangle$$
 and by
$$\gamma_3(s,\tau)= \langle e_m,\int_s^\tau
\{f(\theta,\xi(\theta), \hat y(\theta),\hat z(\theta))-
f(\theta,\xi(\theta),\tilde y(\theta), \tilde
z(\theta))\}d\theta\rangle$$ and rewrite (\ref{4.17}) in the form
$$E[\gamma_1(s,\tau)+\gamma_2(s,\tau)+\gamma_3(s,\tau)]\ge 0.$$

Assume that there exists a number $\delta_0<0$ such that ${\cal
K}^{u,\phi}(s,x)<\delta_0$ and
$$\tau_1=\inf\{\theta\in [s,s+\alpha]:{\cal K}^{y(\theta),z(\theta)}(\theta,\xi(\theta))\le \delta_0\}\wedge \tau.$$
By assumption (\ref{4.17}) holds for  $\tau$ and hence for $\tau_1$.
But this leads to a contradiction since
$$0>\delta_0E(\tau_1-s)\ge E\left
[\int_s^{\tau_1}{\cal
N}^{y(\theta),z(\theta)}(\theta,\xi(\theta))d\theta\right]\ge 0.$$

It remains to check that $\gamma_2(s,s+\Delta s)\to 0$ and
$\gamma_3(s,s+\Delta s)\to 0$ as $\Delta s\to 0$ a.s.

Note that $\gamma_2(s,s+\Delta s)\to 0$ a.s. by  definition of $f$,
properties of $\Gamma(s,t)$
  and uniqueness of a BSDE solution.

Finally we check that $\gamma_3(s,s+\Delta s)\to 0$ as $\Delta s\to
0$ a.s. Note that  the couple $( \tilde y(t), \tilde z(t)), s\le
t\le
 s+\Delta s $ satisfies
\begin{equation}\label{4.12} \tilde y(t)=\Gamma^{*}(s+\Delta s)\phi(s+\Delta s,\xi_{s,x}(s+\Delta s))+
\int_t^{s+\Delta s} f(\theta,\xi_{s,x}(\theta), \tilde
y(\theta),\tilde z(\theta))d\theta-
\end{equation} $$-\int_t^{s+\Delta s} \tilde z(\theta)dw(\theta).$$

Given $ \theta \in[s, s+\Delta s]$,
$${\cal K}^{u,\phi}(\theta,x)=\left(
\frac{\partial\phi}{\partial \theta}+{\cal A}\phi\right)(\theta,
x)+g(\theta,x,\phi(\theta,x), K(u,\nabla \phi)(\theta,x)),$$ set
$$\upsilon(\theta)= \tilde y(s+\Delta s)- \Gamma^{*}(s,\theta)\phi(\theta,\xi_{s,x}(\theta))-
\int_\theta^{s+\Delta s} {\cal
K}^{u,\phi}(\vartheta,\xi_{s,x}(\vartheta))d\vartheta$$ and
$$
\varpi(\theta)= \tilde z(\theta)-K(u,\nabla\phi)
(\theta,\xi_{s,x}(\theta))).$$
 Applying Ito's formula we derive
BSDE to govern the couple $( \upsilon(\theta),\varpi(\theta))$
\begin{equation}\label{4.13}
\upsilon(\theta)=\Gamma^{*}(s,s+\Delta s)\phi(s+\Delta
s,\xi_{s,x}(s+\Delta
s))-\Gamma^{*}(s,\theta)\phi(\theta,\xi_{s,x}(\theta))+\end{equation}
$$+\int_\theta^{s+\Delta s} f(\vartheta,\xi_{s,x}(\vartheta), \tilde y(\vartheta),\tilde z(\vartheta))
d\vartheta - \int_\theta^{s+\Delta s}{\cal
K}^{u,\phi}(\vartheta,\xi_{s,x}(\vartheta))d\vartheta-
$$
$$
-\int_\theta^{s+\Delta s} \tilde
z(\vartheta)dw(\vartheta)+\int_\theta^{s+\Delta s}K(u,\nabla
\phi)(\vartheta,\xi_{s,x}(\vartheta)))dw(\vartheta)=$$
$$\int_\theta^{s+\Delta
s} f(\vartheta,\xi_{s,x}(\vartheta),
\upsilon(\vartheta)+\Gamma^{*}(s,\vartheta)\phi(\vartheta,\xi_{s,x}(\vartheta))+
$$$$+\int_\vartheta^{s+\Delta s} {\cal K}^{u,\phi}(r,\xi_{s,x}(r))dr,\,\varpi(\vartheta)+
K(u, \nabla\phi )(\vartheta,\xi_{s,x}(\vartheta)))d\vartheta+$$
$$
+\int_\theta^{s+\Delta s}\left[\left(\frac{\partial \phi}{\partial
\vartheta}+{\cal A}\phi\right)(\vartheta,\xi_{s,x}(\vartheta))-
{\cal K}^{u,\phi}(\vartheta,\xi_{s,x}(\vartheta))\right]d\vartheta
-\int_\theta^{s+\Delta s} \varpi(\vartheta)dw(\vartheta).$$ We
verify that $(\upsilon,\varpi)$   converges to $(0,0)$   as  $\Delta
s\to 0$. Keeping in mind the estimates for the  generator  $g$   by
standard reasoning based on the Ito formula and the Burkholder
inequality  we can prove that
  $$E\left[\sup_{t\in[s,s+\Delta s]}|\upsilon (t)|^2\right]+
  E\left[\int_s^{s+\Delta s}\|\varpi (\theta)\|^2d\theta\right]\le
   LE\left[\int_s^{s+\Delta s}\|m(\theta,\Delta s)\|^2d\theta\right],$$
where
$$m(\theta,\Delta s)=-{\cal K}^{u,\phi}(\theta,\xi_{s,x}(\theta))+\left(\frac{\partial\phi}{\partial \theta}+
{\cal A}\phi\right)(\theta,\xi_{s,x}(\theta))+$$
$$f(\theta,\xi_{s,x}(\theta),
\upsilon(\theta)+\Gamma^{*}(s,\theta)\phi(\theta,\xi_{s,x}(\theta))+
$$$$+\int_\theta^{s+\Delta s} {\cal K}^{u,\phi}(r,\xi_{s,x}(r))dr,\,\varpi(\theta)+
K(u, \nabla\phi )(\theta,\xi_{s,x}(\vartheta)).$$

Furthermore, since   $\sup_{\theta\in[s,s+\Delta
s]}E[\|\xi_{s,x}(\theta)-x\|^2]\to 0$     as  $\Delta s\to 0$ and
parameters of stochastic equations as well as the function   $\phi$
and its derivatives  are uniformly continuous  in  $x$,  we obtain
$$\lim_{\Delta s \to 0}\sup_{s\le \theta \le s+\Delta s}E[\|m(\theta,\Delta s)\|^2]=0.$$
Hence,
\begin{equation}\label{4.14}
E\left[\sup_{s\le \theta \le s+\Delta s}|\upsilon(\theta)|^2\right]+
E\left[\int_s^{s+\Delta s}\|\varpi (\theta)\|^2d\theta\right]\le \end{equation}
$$ LE\left[\int_s^{s+\Delta s}\|m(\theta,\Delta t)\|^2d\theta\right]\le \varepsilon(\Delta s) \Delta s,$$
 where $\varepsilon(\Delta s)\to 0$   as      $\Delta s\to 0$.  As a result we get that
 $ \tilde y(\theta)$      converges to $\phi(s,x)$  and
 $\tilde z(\theta)$  converges to $[Cu](s,x)+[\nabla\phi A](s,x)$
a.s. as $\Delta s\to 0.$

This estimate does not satisfy yet our purposes.  To get a more
satisfactory estimate  we evaluate the conditional  expectation of
both sides  of  (\ref{4.13}), that leads to
$\upsilon(\theta)=E\left[\int_\theta^{s+\Delta s}n(\vartheta,\Delta
s)d\vartheta|{\cal F}_{\theta}\right],$ where
$$n(\theta,\Delta s)=-{\cal K}^{u,\phi}(\theta,\xi_{s,x}(\theta))+\left[\frac{\partial\phi}{\partial \theta}+
{\cal A}\phi\right](\theta,\xi_{s,x}(\theta))+
f(\theta,\xi_{s,x}(\theta), \tilde y(\theta), \tilde z(\theta))=$$
$$=f(\theta,\xi_{s,x}(\theta), \tilde y(\theta),
 \tilde z(\theta))-$$
$$-f\left(\theta,\xi_{s,x}(\theta),\Gamma^{*}(\theta)\phi(\theta,\xi_{s,x}(\theta))+
\int_\theta^{s+\Delta s}{\cal
K}^{u,\phi}(\theta,\xi_{s,x}(\theta))d\theta, \, K(u,
\nabla\phi)(\theta,\xi_{s,x}(\theta))\right).$$ By Lipschitz
continuity of $f$ we have for $s\le \theta\le s+\Delta s$,
$\|n(\theta,\Delta s)\|\le L[\|\upsilon (\theta)\|$ +
$\|\varpi(\theta)\|],$
 that is  $\|n(\theta,\Delta s)\| \to 0$  a.s. as $ \Delta s\to 0$.

 Hence we have proved that
    $u(s,x)$ is a viscosity subsolution of the Cauchy problem (\ref{4.5}).
    In a similar way we prove that $u(s,x)$ is a supersolution of   (\ref{4.5}) and hence a viscosity solution of this problem.

{\bf Acknowledgement}
Financial support of grant RFBR 12-01-00427-a and  the Minobrnauki
project 1.370.2011
  is gratefully acknowledged

\end{document}